\documentclass[12pt]{amsart}
\usepackage{enumerate}
\usepackage[a4paper,margin=3cm]{geometry}
\usepackage{tikz}
\usepackage{hyperref}
\hypersetup{
    pdftitle=   {Rank-width: Algorithmic and structural results},
   pdfauthor=  {Sang-il Oum}
}

\newcommand\abs[1]{\lvert #1\rvert}
\newtheorem{THM}{Theorem}[section]

\newtheorem{QUE}{Question}
\theoremstyle{remark}

\newcommand\poly{\operatorname{poly}}
\newcommand\pivot\wedge
\newcommand\F{\mathbb F}
\newcommand\rw{\operatorname{rw}}
\newcommand\lrw{\operatorname{lrw}}
\newcommand\cw{\operatorname{cw}}

\newcommand\nlc{\operatorname{nlc}}
\newcommand\boolw{\operatorname{boolw}}
\theoremstyle{definition}

\begin{document}
\title{Rank-width: Algorithmic and structural results}
\author{Sang-il Oum}
\thanks{Supported by Basic Science Research
  Program through the National Research Foundation of Korea (NRF)
  funded by  the Ministry of Science, ICT \& Future Planning
  (2011-0011653).}
\address{Department of Mathematical Sciences, KAIST, Daejeon, South Korea}
\email{sangil@kaist.edu}
\date{\today}
\begin{abstract}
Rank-width is a width parameter of graphs describing whether it is possible to decompose a graph into a tree-like structure by `simple' cuts.
This survey aims to summarize known algorithmic and structural results on rank-width of graphs.
\end{abstract}
\maketitle

\section{Introduction}\label{sec:intro}
Rank-width was introduced by Oum and Seymour \cite{OS2004}. 
Roughly speaking, the rank-width of a graph is the minimum integer $k$ such that the graph can be decomposed into a tree-like structure by recursively splitting its vertex set so that each cut induces a matrix of rank at most $k$. The precise definition will be presented in the next section.

Tree-width is a better known and well-studied width parameter of graphs introduced by Robertson and Seymour~\cite{RS1986a}. 
Tree-width has been very successful for deep theoretical study of graph minor structures
as well as algorithmic applications. In particular, many NP-hard problems can be solved efficiently if the input graph belongs to a class of graphs of bounded tree-width.
However, graphs of bounded tree-width have bounded average degree and therefore the application of tree-width is mostly limited to `sparse' graph classes. 

Clique-width \cite{CER1993,CO2000} aims to extend tree-width by allowing more `dense' graphs to have small clique-width. There are many efficient algorithms based on dynamic programming assuming that a clique-width decomposition is given, see~\cite{HOSG2006}. However, there is no known polynomial-time algorithm to determine whether the clique-width of an input graph is at most $k$ for fixed $k\ge 4$. (There is a polynomial-time algorithm to decide clique-width at most $3$ by Corneil et al.~\cite{CHLRR2012}.)

Rank-width solves this dilemma in some way; there is a polynomial-time algorithm to decide whether the rank-width of an input graph is at most $k$ for each fixed $k$.
In addition, for each graph of rank-width $k$ and clique-width $c$, the inequality \[ 
k\le c\le 2^{k+1}-1\] holds \cite{OS2004} and therefore a class of graphs has bounded clique-width if and only if it has bounded rank-width.
Furthermore a rank-decomposition of width $k$ can be translated into a clique-width decomposition of width at most $2^{k+1}-1$~\cite{OS2004} so that we can use known algorithms for graphs of bounded clique-width.
Thus, if we could extend theorems on tree-width to rank-width, then it becomes useful for some dense graphs.
Moreover there is a correspondence between branch-width of binary matroids and rank-width of bipartite graphs in terms of vertex-minors and pivot-minors, discussed in \cite{Oum2004}. This correspondence allows us to apply the matroid theory to prove statements on rank-width, for instance in \cite{CO2004}.

In 2008, the author with  Hlin\v{e}n\'y, Seese, and Gottlob wrote a survey~\cite{HOSG2006} on various width parameters including rank-width. Our (incomplete) survey aims to supplement that by focusing on algorithmic problems of finding rank-width and structural aspects related to rank-width.

\section{Definitions and equivalent width parameters}
All graphs in this survey are assumed to be simple, meaning that they have no parallel edges and no loops.
We focus on finite undirected graphs in this article,
though Kant\'e and Rao \cite{KR2013} worked on rank-width of directed graphs and edge-colored graphs and Courcelle~\cite{Courcelle2014} extended rank-width and linear rank-width to countable graphs.

In later sections, we discuss and use relations to matroids and branch-width.
For the definition of matroids, please refer to the book of Oxley~\cite{Oxley1992}.
Robertson and Seymour~\cite{RS1991} introduced the notion of branch-width for both graphs and matroids. 
Key concepts of binary matroids related to rank-width can be found in \cite{Oum2004}.

\subsection{Rank-width}
In order to define rank-width, we first need a measure to describe the complexity of a cut, which we call the \emph{cut-rank function} $\rho_G$.
For a graph $G$ and a subset $X$ of the vertex set $V(G)$ of $G$, 
let us define $\rho_G(X)$ to be the rank of a $\abs{X}\times\abs{V(G)-X}$  $0$-$1$ matrix $A_X$ over the binary field where
the entry of $A_X$ on the $i$-th row and $j$-th column is $1$ if and only  if 
the $i$-th vertex in $X$ is adjacent to the $j$-th vertex in $V(G)-X$.
(If $X=\emptyset$ or $X=V(G)$, then $\rho_G(X)=0$.)
Note that this matrix is over the binary field; thus, for instance, the rank of $\left(\begin{smallmatrix}1&1&0\\
0&1&1\\
1&0&1\end{smallmatrix} \right)$ is equal to $2$, as $1+1=0$.
The cut-rank function satisfies the following inequalities.
\begin{align}
  \rho_G(X)&=\rho_G(V(G)-X) \text{ for all }X\subseteq V(G),
             \label{eq:symmetric}
  \\
  \rho_G(X)+\rho_G(Y)&\ge \rho_G(X\cup Y)+\rho_G(X\cap Y) \text{ for all }X,Y\subseteq V(G).
                       \label{eq:submodular}
\end{align}
The inequality \eqref{eq:submodular} is called the \emph{submodular} inequality.

\medskip
A tree is \emph{subcubic} if every node has degree $1$ or $3$.
A \emph{rank-decomposition} of a graph $G$ is a pair $(T,L)$ of a subcubic tree $T$ with at least two nodes and a bijection $L$ from $V(G)$ to  the set of all leaves of $T$.
For each edge $e$ of $T$, $T-e$ induces a partition $(A_e,B_e)$ of the leaves of $T$ and we say that the \emph{width} of $e$ is $\rho_G(L^{-1}(A_e))$. 
By \eqref{eq:symmetric}, the choice of $A_e$, $B_e$ does not change the width of $e$.
The \emph{width} of a rank-decomposition $(T,L)$ is the maximum width of edges in $T$.

The \emph{rank-width} of a graph $G$, denoted by $\rw(G)$, is the minimum width over all rank-decompositions of $G$.  
(If $G$ has less than two vertices, then $G$ has no rank-decompositions; in this case we say that $G$ has rank-width $0$.)

\subsection{Linear rank-width}
Let us define a linearized variant of rank-width as follows. A \emph{linear layout} of an $n$-vertex graph $G$ is a permutation $v_1,v_2,\ldots,v_n$ of the vertices of $G$. The \emph{width} of a linear layout is $\min_{i=1}^{n-1}\rho_G(\{v_1,v_2,\ldots,v_i\})$. The \emph{linear rank-width} of a graph $G$, denoted by $\lrw(G)$, is the minimum width over all linear layouts of $G$. (If $G$ has less than two vertices, then we say the linear rank-width of $G$ is $0$.)

An alternative way to define the linear rank-width is to define a linear rank-decomposition $(T,L)$ to be a rank-decomposition such that $T$ is a caterpillar 
and then define linear rank-width as the minimum width over all linear rank-decompositions.   (A \emph{caterpillar} is a tree having a path such that every leaf of the tree is either in the path or adjacent to a vertex in the path.)
This immediately implies that the rank-width is always less than or equal to the linear rank-width.

\subsection{\emph{Equivalent} width parameters}
Let us write $\cw(G)$ to denote the clique-width, introduced in~\cite{CER1993,CO2000}. 
As we recall in Section~\ref{sec:intro}, a class of graphs has bounded rank-width if and only if it has bounded clique-width by the following inequality. 
\begin{THM}[Oum and Seymour~\cite{OS2004}]\label{thm:rwcw}
$  \rw(G)\le \cw(G)\le 2^{\rw(G)+1}-1$.
\end{THM}
One can also prove, by the same idea, that 
linear clique-width of a graph $G$ is at least $\lrw(G)$ and at most $2^{\lrw(G)}+1$.

The inequalities in Theorem~\ref{thm:rwcw} are essentially tight:
\begin{itemize}
\item The $n\times n$ grid has rank-width $n-1$, shown by Jel{\'{\i}}nek~\cite{Jelinek2010} 
and clique-width $n+1$, shown by Golumbic and Rotics~\cite{GR2000}.
\item 
For every $k$, there exists a graph $G$ such that its rank-width is at most $k+1$ and the clique-width is at least $2^{\lfloor k/2\rfloor-1}$. This is due to a theorem by Corneil and Rotics~\cite{CR2005}, that constructs a graph of tree-width $k$ whose clique-width is at least $2^{\lfloor k/2\rfloor -1}$. Oum~\cite{Oum2006c} showed that a graph of tree-width $k$ has rank-width at most $k+1$.
\end{itemize}

We may say that clique-width  and rank-width are \emph{equivalent} in the sense that one is bounded if and only if the other is bounded.
There are a few other width parameters of graphs that are equivalent to rank-width.
We list them below.
\begin{itemize}
\item NLC-width was introduced by Wanke~\cite{Wanke1994} and is convenient to design an algorithm based on dynamic programming on graphs of bounded clique-width or rank-width. Let us write $\nlc(G)$ to denote the NLC-width of $G$.
Johansson~\cite{Johansson1998,Johansson2001b} proved that \[\cw(G)\le \nlc(G)\le 2\cw(G).
\]
It can be easily shown that \[
\rw(G)\le \nlc(G)\le 2^{\rw(G)}.
\]

\item Boolean-width, denoted by $\boolw(G)$, was introduced by Bui-Xuan, Telle, and Vatshelle~\cite{BTV2011}. They showed that 
\begin{align*}
\log_2 \cw(G) &\le \boolw(G)\le \cw(G) \text{ and }\\
\log_2 \rw(G)&\le \boolw(G)\le \frac{\rw(G)^2}{4}+O(\rw(G)).
\end{align*}
In addition, they showed that for each $k$, there exist graphs $G$ and $H$ such that the 
$\rw(G)\ge k$ and $\boolw(G)\le 2\log_2 k+4$
and $\rw(H)=k$ and $\boolw(H)\ge (k+1)/6$. It is not yet known~\cite{BTV2011} whether $\rw(G)^2$ can be reduced to $O(\rw(G))$ in the above inequality in general.

The definition of boolean-width is a slight variation of rank-width, only changing the cut-rank function into an alternative function called the \emph{boolean-dimension}. For a graph $G$ and  a subset $A$ of vertices of $G$, let $U(A)$ be the set of all subsets of $V(G)-A$ that can be represented as $N(X)-A$ for some $X\subseteq A$, where $N(X)$ is the set of all vertices in $V(G)-X$ that are adjacent to a vertex in $X$. Then the boolean-dimension of $A$ is defined to be $\log_2\abs{ U(A)}$. It turns out that $\abs{U(A)}=\abs{U(V(G)-A)}$ (see \cite{Kim1982}) and therefore the boolean-dimension is symmetric. 
However, it is not necessarily integer-valued.

\item One may consider defining rank-width over a different field other than the binary field $GF(2)$ by changing the definition of cut-rank function.
For a field $\F$, let us say that $\F$-rank-width of a graph $G$, denoted by $\rw_\F(G)$, is the obtained width parameter from rank-width by changing the cut-rank function to evaluate the rank of a $0$-$1$ matrix over $\F$ instead of the binary field. 
The proof of Oum and Seymour~\cite{OS2004} can be extended easily to $\F$-rank-width to prove
\[
\rw_\F(G)\le \cw(G)\le 2^{\rw_\F(G)+1}-1.
\]
Thus rank-width and $\F$-rank-width are equivalent.

Furthermore 
Oum, S\ae ther, and Vatshelle~\cite{OSV2014} observed that 
if $\F$ is the field  of all rational numbers
(or more generally, $\F$ is a field of characteristic $0$ or $2$),
then one can enhance the above inequality by showing
\[
\rw(G)\le \rw_{\mathbb F}(G)\le \cw(G)\le 2^{\rw(G)+1}-1.
\]

\end{itemize}
Most of other width parameters mentioned above
do not have a direct algorithm to find a corresponding decomposition and therefore we simply employ algorithms on finding rank-decompositions and convert them.
But $\mathbb F$-rank-width is similar to rank-width in the sense that 
one can use the theorem of Oum and Seymour~\cite{OS2004} to design its own fixed-parameter approximation algorithm.
This was used to design a faster algorithm for certain vertex partitioning problems \cite{OSV2014}. 
When Kant\'e and Rao~\cite{KR2013} extended rank-width to directed graphs and edge-colored graphs, they use the idea of $\F$-rank-width.

\section{Discovering Rank-width and linear rank-width}\label{sec:discover}
We skip the discussion on algorithms on graphs of small rank-width, 
except the important meta theorem of Courcelle, Makowsky, and Rotics~\cite{CMR2000}. 
Their meta theorem states the following.
\begin{THM}[Courcelle, Makowsky, and Rotics~\cite{CMR2000}]
  Let $k$ be a fixed constant.
  For every closed monadic second-order formula of the first kind on graphs,
  there is a $O(n^3)$-time algorithm to determine whether the input graph
  of rank-width at most $k$ satisfies the formula.
\end{THM}
Here, the \emph{monadic second-order formula of the first kind} on graphs 
are logic formulas which allow $\exists$, $\forall$, $\lnot$, $\land$, $\lor$, $\in$, $($, $)$, 
$\mathbf{true}$, and $\mathbf{adj}$ with first-order variables each representing a vertex and 
set variables each representing a set of vertices, where $\mathbf{adj}(x,y)$ is true if two vertices $x$ and $y$ are adjacent.
If such a formula has no free variable, then it is called \emph{closed}.
We refer the readers to the previous survey \cite{HOSG2006}.
We would like to mention that Grohe and Schweitzer~\cite{GS2015} announced a polynomial-time algorithm to check whether two graphs are isomorphic when the input graphs have bounded rank-width.

We often talk about \emph{fixed-parameter tractable (FPT)} algorithms. This means that for a parameter $k$, the running time of the algorithm is at most $f(k)n^c$ for some $c$ where $n$ is the length of the input. For more about parametrized complexity, please see the following books \cite{DF2013,CFKLMPPS2015}. For us, the parameter is usually the rank-width of the input graph. A fixed-parameter tractable algorithm is thought to be better than XP algorithms, that is an algorithm whose running time is bounded above by $f(k)n^{g(k)}$. 

\subsection{Hardness}
Computing rank-width is NP-hard, as mentioned in \cite{Oum2006}. This can be deduced easily from the following three known facts.

\begin{enumerate}
\item Oum~\cite{Oum2004} showed that the branch-width of a binary matroid is equal to the rank-width of its fundamental graph plus one. Note that every cycle matroid is binary.
\item Hicks and McMurray Jr.~\cite{HM2005} and independently Mazoit and Thomass\'e~\cite{MT2005} showed that if a graph $G$ is not a forest, then the branch-width of  the cycle matroid $M(G)$ is equal to the branch-width of $G$.
\item Seymour and Thomas \cite{ST1994} proved that computing branch-width of graphs is NP-hard. 
\end{enumerate}

It is also NP-hard to compute linear rank-width. This is implied by the following facts.
\begin{enumerate}
\item Kashyap~\cite{Kashyap2008} proved that it is NP-hard to compute the path-width of a binary matroid given by its representation. 
\item By Oum~\cite{Oum2004}, it is straightforward to deduce that the path-width of a binary matroid is equal to the linear rank-width of its fundamental graph plus one. 
\end{enumerate}

\subsection{Computing rank-width}

Oum and Seymour~\cite{OS2005} presented a generic XP algorithm for 
symmetric submodular integer-valued functions. Their algorithm implies that in time $O(n^{8k+12}\log n)$ we can find a rank-decomposition of width at most $k$, if it exists, for an input $n$-vertex graph.

But we definitely prefer to remove $k$ in the exponent of $n$.
Courcelle and Oum~\cite{CO2004} constructed a fixed-parameter tractable algorithm to decide, for fixed $k$, whether the rank-width is at most $k$, in time $O(g(k)n^3)$. However their algorithm does not provide a rank-decomposition of width at most $k$ directly, because their algorithm uses forbidden vertex-minors, which we discuss in the next section.

This problem was solved later. Hlin\v en\'y and Oum~\cite{HO2006} devised an algorithm to find a rank-decomposition of width at most $k$ if it exists, for an input $n$-vertex graph in time $O(g(k)n^3)$. Here $g(k)$ is a huge function. The proof uses the correspondence between branch-width of binary matroids and rank-width of graphs and still uses a huge list of forbidden minors in matroids.

\medskip

For bipartite circle graphs, one can compute rank-width in polynomial time. This is due to the following results.
\begin{itemize}
\item Seymour and Thomas~\cite{ST1994} proved that branch-width of planar graphs can be computed in polynomial time.
\item Hicks and McMurray Jr.~\cite{HM2005} and independently Mazoit and Thomass\'e~\cite{MT2005} showed that if a graph $G$ is not a forest, then the branch-width of  the cycle matroid $M(G)$ is equal to the branch-width of $G$.
\item Oum~\cite{Oum2004} showed that the branch-width of a binary matroid is equal to the rank-width of its fundamental graph plus one. Note that every cycle matroid is binary.
\item De Fraysseix~\cite{deFraysseix1981} proved that a graph is a bipartite circle graph if and only if it is a fundamental graph of the cycle matroid of a planar graph.
\end{itemize}
It would be interesting to extend the class of graphs on which rank-width can be computed exactly.
\begin{QUE}
  Can we compute rank-width of circle graphs in polynomial time?
\end{QUE}

\subsection{Fixed-parameter approximations}
For many algorithmic applications, we restrict our input graph to have bounded rank-width.
Most of such algorithms assume that the input graph is given with a rank-decomposition of small width and for such algorithms, it is usually enough to have an efficient algorithm that outputs a rank-decomposition of width at most $f(k)$ for  some function $f$
or confirms that the rank-width of the input graph is larger than $k$.

The first such algorithm was devised by Oum and Seymour~\cite{OS2004} with $f(k)=3k+1$ and the running time $O(8^k n^9\log n)$, where $n$ is the number of the vertices in the input graph.
It was later improved by Oum~\cite{Oum2006} with $f(k)=3k+1$ and the running time $O(8^kn^4)$.
This allows us to compute rank-width within a constant factor in polynomial time when the input graph is known to have rank-width at most $c\log n$ for some constant $c$.

In the same paper, Oum~\cite{Oum2006} presented another algorithm with $f(k)=3k-1$ whose running time is $O(g(k) \,n^3)$ for some huge function $g(k)$.
 \begin{QUE}
   Does there exist an algorithm with a function $f(k)$ and a constant $c$ that finds a rank-decomposition of width at most $f(k)$ or confirms that the rank-width of the $n$-vertex input graph is larger than $k$,
 in time $O(c^k n^3)$?
 \end{QUE}

One may hope to find a faster approximation algorithm for rank-width. 
For tree-width, 
Bodlaender, Drange, Dregi, Fomin, Lokshtanov, and Pilipczuk~\cite{BDDFLP2016}
presented an algorithm that, for an $n$-vertex input graph and an integer $k$, finds
a tree-decomposition of width at most $5k+4$ or confirms that the tree-width is larger than $k$ in time $2^{O(k)}n$.
However, for rank-width, it seems quite non-trivial to reduce $n^3$ into $n^c$ for some $c<3$. This may be justified by the observation that one needs to compute the rank of matrices and computing a rank of a low-rank matrix may be one of the bottlenecks. There are some algorithms computing rank faster \cite{CKL2013} but it is not clear whether those algorithms can be applied to answer the following question.
 \begin{QUE}
   Does there exist an algorithm with functions $f(k)$, $g(k)$ and a constant $c<3$ that finds a rank-decomposition of width at most $f(k)$ or confirms that the rank-width of the $n$-vertex input graph is larger than $k$,
 in time $O(g(k)n^c)$? 
 \end{QUE}

\subsection{Computing linear rank-width}
Let us turn our attention to linear rank-width.
Nagamochi~\cite{Nagamochi2012} presented a simple generic XP algorithm analogous to Oum and Seymour~\cite{OS2005} for path-width. His algorithm implies that we can find a linear rank-decomposition of width at most $k$ by calling the algorithm to minimize submodular functions  
 \cite{Orlin2009,Iwata2003,IFF2001,Schrijver2000} at most  $O(n^{2k+4})$ times, where $n$ is the number of vertices of the input graph. In fact, for the linear rank-width, one can use a direct combinatorial algorithm in \cite{Oum2006} to avoid using the generic algorithm to minimize submodular functions and improve its running time.

Until recently there was no known algorithm that, in time $O(g(k)n^c)$, finds a linear rank-decomposition of width at most $k$ or confirms that a graph has linear rank-width larger than $k$ for an input $n$-vertex graph.
By the well-quasi-ordering property to be discussed in Subsection~\ref{subsec:wqo}, there exists an algorithm that decides in time $O(g(k)n^c)$ whether the input $n$-vertex graph has linear rank-width at most $k$, but such an algorithm would not easily produce a linear rank-decomposition of small width.
Hlin\v en\'y~\cite{Hlineny2016} claimed that one can use the  self-reduction technique \cite{FL1988} to convert a decision algorithm for linear rank-width at most $k$ into a construction algorithm.

Recently, Jeong, Kim,  and Oum \cite{JKO2016} devised a direct algorithm 
for graphs of rank-width at most $k$
that runs in time $O(g(k)^\ell n^3)$ 
to output a linear rank-decomposition of width at most $\ell$ 
or confirm that the linear rank-width is larger than $\ell$.
Now observe the following inequality shown by Kwon~\cite{Kwon2015} (also in Adler, Kant\'e, and Kwon \cite{AKK2015}).
\begin{THM}[Kwon~\cite{Kwon2015}; Adler, Kant\'e, and Kwon \cite{AKK2015}]
Let $G$ be an $n$-vertex graph. 
Then $\lrw(G)\le \rw(G)\lfloor \log_2 n\rfloor$. 
\end{THM}
Thus $\ell\le k\lfloor \log_2 n\rfloor$ and so we can conclude the following.
\begin{THM}[Jeong, Kim, and Oum~\cite{JKO2016}]\label{thm:lrwalgo}
  For each fixed integer $k$, there is a polynomial-time (XP) algorithm to compute the linear rank-width of an input graph of rank-width at most $k$.
\end{THM}
However, it is not known whether there is a fixed-parameter tractable algorithm to compute linear rank-width of a graph of rank-width $k$. 

Theorem~\ref{thm:lrwalgo} generalizes the polynomial-time algorithm of Adler, Kant\'e, and Kwon~\cite{AKK2015} that computes the linear rank-width of graphs of rank-width at most~$1$.

\subsection{Exact algorithms}
Though the number of rank-decompositions of an $n$-vertex graph is a superexponential function, we can still determine the rank-width of an $n$-vertex graph in time $O(2^n \poly(n))$ where $\poly(n)$ is some polynomial in $n$, shown by Oum~\cite{Oum2009}.  This algorithm can find an optimal rank-decomposition as well.
\begin{QUE}
  Is it possible to compute rank-width of an $n$-vertex graph in time $O(c^n)$ for some $c<2$?
\end{QUE}
It is easy to compute the linear rank-width of an $n$-vertex graph in time $O(2^n\poly(n))$
by dynamic programming \cite{BFKKT2012}.
All algorithms mentioned so far in this subsection use exponential space. 

\subsection{Software}

Krause implemented a simple dynamic programming algorithm
to compute the rank-width of a graph.\footnote{\url{http://pholia.tdi.informatik.uni-frankfurt.de/~philipp/software/rw.shtml}}
This software is now included in the open source mathematics software
package called SAGE; see the
manual.\footnote{\url{http://www.sagemath.org/doc/reference/graphs/sage/graphs/graph_decompositions/rankwidth.html}}

Friedmansk\'y wrote the Master thesis \cite{Friedmansky2011}
about implementing the exact exponential-time algorithm of Oum~\cite{Oum2009}.

Bui-Xuan, Raymond, and Tr\'ebuchet~\cite{BRT2011} implemented, in SAGE, the algorithm of Oum~\cite{Oum2006} that either outputs a rank-decomposition of width at most $3k+1$ or confirms that rank-width is larger than $k$ in time $O(n^4)$ for an input $n$-vertex graph.

\subsection{Other results}

Lee, Lee, and Oum \cite{LLO2010} showed that asymptotically almost surely  the Erd\H os-R\'enyi random graph $G(n,p)$ has rank-width $\lceil n/3\rceil -O(1)$
if $p$ is a constant between $0$ and $1$.
Furthermore, if  $\frac{1}{n}\ll p \le \frac{1}{2}$, then the rank-width is $\lceil \frac{n}{3}\rceil-o(n)$, 
and if $p = c/n$ and $c > 1$, then  rank-width is at least $r n$ for some $r = r(c)$.

\section{Structural aspects}\label{sec:structure}
As tree-width is closely related to graph minors, 
rank-width is related to vertex-minors and pivot-minors, observed in~\cite{Oum2004}.
 We first review the definition of pivot-minors and vertex-minors
and then discuss various theorems, mostly motivated by the theory of tree-width and graph minors.

\subsection{Vertex-minors and pivot-minors}
For a graph $G$ and its vertex $v$,  $G*v$  is the graph such that $V(G*v)=V(G)$ and two  vertices $x$, $y$ are adjacent in $G*v$ if and only if 
\begin{enumerate}[(i)]
\item 
both $x$ and $y$ are neighbors of $v$ in $G$ and $x$ is non-adjacent to $y$, or 
\item $x$ is adjacent to $y$ and $v$ is non-adjacent to at least one of $x$ and $y$.
\end{enumerate}
Such an operation is called the \emph{local complementation} at $v$. 
Clearly $G*v*v=G$.
Two graphs are \emph{locally equivalent} if one is obtained from another by a sequence of local complementations.
We say that a graph $H$ is a \emph{vertex-minor} of a graph $G$ if $H$ is an induced subgraph of a graph locally equivalent to $G$.

Local complementation is a useful tool to study rank-width of graphs, because local complementation preserves the cut-rank function, thus preserving rank-width and linear rank-width \cite{Oum2004}. This implies that if $H$ is a vertex-minor of $G$, then the rank-width of $H$ is less than or equal to the rank-width of $G$
and the linear rank-width of $H$ is less than or equal to the linear rank-width of $G$.

Pivot-minors are restricted version of vertex-minors.
For an edge $uv$, we define $G\pivot uv=G*u*v*u$. This operation is called a \emph{pivot}. Two graphs are \emph{pivot equivalent} if one is obtained from another by a sequence of pivots. A graph $H$ is a \emph{pivot-minor} of a graph $G$ if $H$ is an induced subgraph of a graph pivot equivalent to $G$. 
There is a close relation between the pivot operation and the operation of switching a base in a binary matroid, which motivates the study of pivot-minors. 
Every pivot-minor of a graph $G$ is also a vertex-minor of $G$ but not vice versa.

Vertex-minors of graphs can be seen as minors of isotropic systems, a concept introduced by Bouchet~\cite{Bouchet1988,Bouchet1989a,Bouchet1987a}. Similarly pivot-minors of graphs can be seen as minors of binary delta-matroids, also introduced by Bouchet~\cite{Bouchet1987d,Bouchet1988b}.

\subsection{Relation to tree-width and branch-width}
Let $G$ be a graph whose tree-width is $k$.
Kant\'e~\cite{Kante2006a} showed that rank-width  of $G$ is at most $4k+2$ by showing how to convert a tree-decomposition into a rank-decomposition constructively.
This inequality was improved by Oum~\cite{Oum2006c}, who showed that rank-width is at most $k+1$.

In fact, the following are shown.
The \emph{incidence graph} $I(G)$ of a graph $G$ is the graph obtained from $G$ by subdividing each edge exactly once.
\begin{THM}[Oum~\cite{Oum2006c}]
  Let $G$ be a graph of branch-width $k$.  Then the incidence graph $I(G)$ of $G$ has rank-width $k$ or $k-1$, unless $k=0$.
\end{THM}
The converse does not hold; the complete graph $K_n$ has tree-width $n-1$ and rank-width $1$ if $n\ge 2$.
We remark that Courcelle~\cite{Courcelle2015} proved that if $G$ is a graph of tree-width $k$, then its incidence graph $I(G)$ has clique-width at most $2k+4$.

Since $G$ is a vertex-minor of $I(G)$, we deduce that the rank-width of $G$ is less than or equal to the branch-width of $G$, unless $G$ has branch-width $0$. This inequality can be generalized to an arbitrary field.
\begin{THM}[Oum, S{\ae}ther, Vatshelle~\cite{OSV2014}]
  For every field $\F$ and a graph $G$, $\rw_\F(G)$ is less than or equal to branch-width of $G$, unless $G$ has branch-width $0$.
\end{THM}

The line graph $L(G)$ of a simple graph $G$ is a vertex-minor  of the incidence graph $I(G)$, obtained by applying local complementations at vertices of $G$ in $I(G)$. This implies that the rank-width of $L(G)$ is less than or equal to the rank-width of $I(G)$. 
The following theorem proves that they are very close to each other.
\begin{THM}[Oum~\cite{Oum2006a}]
  Let $G$ be a simple graph of branch-width $k$.  Then the line graph of $G$ has rank-width $k$, $k-1$,  or $k-2$.
\end{THM}
Surprisingly, Adler and Kant\'e~\cite{AK2015} showed that for trees, linear rank-width is equal to path-width.

Another relation to tree-width was discovered by Kwon and Oum~\cite{KO2013}. They have shown that graphs of small rank-width are precisely pivot-minors of graphs of small tree-width.
\begin{THM}[Kwon and Oum~\cite{KO2013}]\label{thm:rwtwpivot}
  \begin{enumerate}
  \item Every graph of rank-width  $k$
is a pivot-minor of a graph of tree-width at most $2k$.
\item Every graph of linear rank-width $k$ is a pivot-minor of a graph of path-width at most $k+1$.
  \end{enumerate}
\end{THM}
Hlin{\v{e}}n\'y, Kwon, Obdr{\v{z}}\'alek, and Ordyniak~\cite{HKOO2016} found an analogue of Theorem~\ref{thm:rwtwpivot} to relate tree-depth and shrub-depth of graphs by vertex-minors.

Fomin, Oum, and Thilikos \cite{FOT2010} showed that when  graphs are planar, or $H$-minor-free, then 
having bounded tree-width is equivalent to having bounded rank-width.
This is already proven in Courcelle and Olariu \cite{CO2000} without explicit bounds because they use logical tools.
For instance, in \cite{FOT2010}, it is shown that 
the tree-width of a planar graph $G$ is at most $72\rw(G)-2$, 
and if a graph $G$
has no $K_r$ minor with $r>2$, then its tree-width is at most $2^{O(r\log\log r)}\rw(G)$.
The last bound  on tree-width can be improved to $2^{O(r)}\rw(G)$
by using a recent result of Lee and Oum~\cite{LO2014}.

\subsection{Chromatic number}
A class $\mathcal I$ of graphs is \emph{$\chi$-bounded} if there exists a function $f$ such that $\chi(G)\le f(\omega(G))$ for all graphs $G\in\mathcal I$. 
 Here, $\omega(G)$ is the maximum size of a clique in $G$ and $\chi(G)$ is the chromatic number of $G$.
Dvo{\v{r}}{\'a}k and Kr{\'a}l'~\cite{DK2012} proved that the class of graphs of rank-width at most $k$ is \emph{$\chi$-bounded} as follows.
\begin{THM}[Dvo{\v{r}}{\'a}k and Kr{\'a}l'~\cite{DK2012}]
  There exists a function $f:\mathbb Z\times\mathbb Z\to \mathbb Z$ such that
for every graph $G$ of rank-width at most $k$,
$\chi(G)\le f(\omega(G),k)$.
\end{THM}
More generally, Geelen (see \cite{DK2012,CKO2015}) conjectured that for each fixed graph $H$, the class of graphs with no $H$ vertex-minor is $\chi$-bounded. This conjecture has been verified when $H$ is the wheel graph on $6$ vertices~\cite{DK2012} and a fan graph~\cite{CKO2015}.
  \subsection{Duality}
As rank-width is an instance of branch-width of some symmetric submodular function, we can use a theory developed for symmetric submodular functions by \cite{RS1991,OS2005,AMNT2009,LMT2010,DO2014}. In particular, one can use the concept called a tangle to verify that a graph has large rank-width.

For a graph $G$, a \emph{$\rho_G$-tangle} of order $k$ is a set $\mathcal T$ of subsets of $V(G)$ satisfying the following three axioms.
\begin{enumerate}[(T1)]
\item For $A\subseteq V(G)$, if $\rho_G(A)< k$, then either $A\in \mathcal T$ or $V(G)-A\in\mathcal T$.
\item If $A,B,C\in \mathcal T$, then $A\cup B\cup C\neq V(G)$.
\item For all $v\in V(G)$, $V(G)-\{v\}\notin \mathcal T$.
\end{enumerate}
A theorem of Robertson and Seymour~\cite{RS1991} implies the following.
\begin{THM}[Robertson and Seymour~\cite{RS1991}]
  For an integer $k$, 
 a graph $G$ has a $\rho_G$-tangle of order $k$ if and only if its rank-width is at least $k$.
\end{THM}

A similar concept can be defined for linear rank-width. 
A \emph{$\rho_G$-obstacle} of  order $k$ is a set $\mathcal O$ of subsets of $V(G)$ satisfying the following three axioms.
\begin{enumerate}[(B1)]
\item For all $A\in \mathcal O$, $\rho_G(A)<k$.
\item If $X\subseteq Y\in \mathcal O$ and $\rho_G(X)<k$, then $X\in \mathcal O$.
\item If  $A\cup B\cup C=V(G)$, $A\cap B=\emptyset$, $\rho_G(A)<k$, $\rho_G(B)<k$ and $\abs{C}\le 1$,
then exactly one of $A$ and $B$ is in $\mathcal O$.
\end{enumerate}
The following theorem is implied by Fomin and Thilikos~\cite[Theorem 3]{FT2003}.
\begin{THM}[Fomin and Thilikos~\cite{FT2003}]
For an integer $k$, 
a graph $G$ has  a $\rho_G$-obstacle of order $k$ if and only if its linear rank-width is at least $k$.
\end{THM}

\subsection{Large rank-width}
Robertson and Seymour~\cite{RS1986} proved that every graph with sufficiently large tree-width contains a large grid graph as a minor. An analogous conjecture was made in~\cite{Oum2006a} as follows.
\begin{QUE}
  Is it true that for each fixed bipartite circle graph $H$,  every graph $G$ with sufficiently large rank-width contains a pivot-minor isomorphic to $H$?
\end{QUE}
This has been confirmed when $G$ is a line graph, a bipartite graph, or a circle graph, see \cite{Oum2006a}.

Adler, Kant\'e and Kwon~\cite{AKK2015a} proved that for each tree $T$, every distance-hereditary graph with sufficiently large linear rank-width contains a vertex-minor isomorphic to $T$. It would be interesting to determine whether every graph of sufficiently large linear rank-width contains a vertex-minor isomorphic to a fixed tree $T$.

\subsection{Well-quasi-ordering}\label{subsec:wqo}
A pair $(Q,\preceq)$ of a set $Q$ and a relation $\preceq$ on $Q$ is a \emph{quasi-order} if 
\begin{enumerate}
\item 
$x\preceq x$ for all $x\in Q$ and
\item $x\preceq z$ if $x\preceq y$ and $y\preceq z$.
\end{enumerate}
A quasi-order $(Q,\preceq)$ is a \emph{well-quasi-order} if every infinite sequence $q_1,q_2,\ldots$ of $Q$ has a pair $q_i$, $q_j$ ($i<j$) such that $q_i\preceq q_j$.
We say that $Q$ is \emph{well-quasi-ordered} by the relation $\preceq$.

Motivated by the celebrated graph minor theorem by Robertson and Seymour~\cite{RS2004a} and its special case on tree-width~\cite{RS1990},
Oum~\cite{Oum2004a} proved that graphs of bounded rank-width are \emph{well-quasi-ordered} under taking pivot-minors as follows.
\begin{THM}[Oum~\cite{Oum2004a}]\label{thm:wqo}
For all positive integers $k$, every infinite sequence $G_1,G_2,\ldots$ of graphs of rank-width at most $k$ admits a pair $G_i$, $G_j$ ($i<j$) such that $G_i$ is isomorphic to a pivot-minor of $G_j$.
\end{THM}
Theorem~\ref{thm:wqo} has been generalized to 
skew-symmetric or symmetric matrices over a fixed finite field (Oum~\cite{Oum2006b})
and
$\sigma$-symmetric matrices over  a fixed finite field (Kant\'e~\cite{Kante2012}).

Courcelle and Oum~\cite{CO2004} constructed a modulo-$2$ counting monadic second-order logic formula to decide whether a fixed graph $H$ is isomorphic to a vertex-minor of an input graph. This can be combined with the meta-theorem of Courcelle, Makowsky, and Rotics~\cite{CMR2000} to deduce the following.
\begin{THM}[Courcelle and Oum~\cite{CO2004}]
  Let $\mathcal I$ be a set of graphs closed under taking vertex-minors.
  If $\mathcal I$ has bounded rank-width, then there exists an algorithm that decides whether an  input $n$-vertex graph is in $\mathcal I$ in time $O(n^3)$.
\end{THM}

We remark that a positive answer to the following question would imply the graph minor theorem of Robertson and Seymour.
\begin{QUE}
  Are graphs well-quasi-ordered under taking pivot-minors?
  In other words, does every infinite sequence $G_1,G_2,\ldots$ of graphs admit a pair $G_i$, $G_j$ ($i<j$) such that $G_i$ is isomorphic to a pivot-minor of $G_j?$
\end{QUE}
We would also like to have an algorithm to detect pivot-minors and vertex-minors of a graph.
\begin{QUE}
  For a fixed graph $H$, can we decide in polynomial time whether an input graph has a pivot-minor (or a vertex-minor) isomorphic to $H$?
\end{QUE}

\subsection{Forbidden vertex-minors}
By Theorem~\ref{thm:wqo},
we can easily deduce that there exists a finite set $F_k$ of graphs such that 
a graph $G$ has rank-width at most $k$ if and only if $F_k$ contains no vertex-minor of $G$, 
because otherwise there will be an infinite anti-chain of  vertex-minor-minimal graphs having rank-width $k+1$. However, the well-quasi-ordering theorem does not provide an upper bound on $\abs{F_k}$.

Oum \cite{Oum2004} showed that pivot-minor-minimal graphs of rank-width $k+1$ have at most $(6^{k+1}-1)/5$ vertices. This gives an explicit upper bound on $\abs{F_k}$ as well as a constructive enumeration algorithm for graphs in $F_k$.

For $k=1$, $F_1$ is known.
Oum~\cite{Oum2004} showed that graphs of rank-width at most $1$ are precisely distance-hereditary graphs; a graph $G$ is \emph{distance-hereditary} if in every connected induced subgraph, a shortest path between two vertices is also a shortest path in $G$ \cite{BM1986}.
(This also proves that distance-hereditary graphs have clique-width at most $3$ as an easy consequence of Theorem~\ref{thm:rwcw}, which was initially proved by Golumbic and Rotics~\cite{GR2000}.)
A theorem of Bouchet~\cite{Bouchet1987a,Bouchet1988a} implies that a graph $G$ has rank-width at most $1$ if and only if $C_5$ (the cycle graph of length $5$) is not isomorphic to a vertex-minor of $G$.
A similar argument can be used to show that a graph $G$ has rank-width at most $1$ if and only if neither $C_5$ nor $C_6$ is isomorphic to a pivot-minor of $G$.

\medskip
For the linear rank-width, we can still use Theorem~\ref{thm:wqo} to deduce that there exists a finite set $L_k$ of graphs such that a graph $G$ has linear rank-width at most $k$ if and only if $L_k$ contains no vertex-minor of $G$.
Unlike the situation in rank-width, so far no upper bound on the size of graphs in $L_k$ is known. 
Adler, Farley, and Proskurowski~\cite{AFP2013} determined $L_1$ as a set of three particular graphs.
Jeong, Kwon, and Oum \cite{JKO2014} proved that $\abs{L_k}\ge 3^{\Omega(3^k)}$.
Adler, Kant\'e and Kwon \cite{AKK2015a} presented an algorithm to construct, for each $k$, a finite set $L'_k$ of distance-hereditary graphs of linear rank-width $k+1$ such that 
$\abs{L_k'}\le 3^{O(3^k)}$, 
every proper vertex-minor of a graph in $L_k'$ has linear rank-width at most $k$, 
and
every distance-hereditary graph of linear rank-width $k+1$ contains at least one graph in $L_k'$.

\subsection*{Acknowledgments}
The author would like to thank Mamadou M. Kant\'e, Eun Jung Kim,  O-joung Kwon,  and Jisu Jeong for their helpful comments on an early draft of this paper.
The author would also like to thank anonymous referees for their helpful suggestions.

\end{document}